%% file: main_revised_3rd_round.tex
\documentclass[pdflatex,sn-mathphys-num]{sn-jnl}% Math and Physical Sciences Numbered Reference Style
\usepackage{graphicx}
\usepackage{multirow}
\usepackage{amsmath,amssymb,amsfonts}
\usepackage{amsthm}
\usepackage{mathrsfs}
\usepackage[title]{appendix}
\usepackage{xcolor}
\usepackage{textcomp}
\usepackage{manyfoot}
\usepackage{booktabs}
\usepackage{algorithm}
\usepackage{algorithmicx}
\usepackage{algpseudocode}
\usepackage{listings}
\usepackage{comment}
\usepackage{hyperref}
\usepackage{subcaption}
\usepackage{caption}
 
\theoremstyle{thmstyleone}
\newtheorem{theorem}{Theorem}%  meant for continuous numbers
\newtheorem{remark}{Remark}

\newtheorem{lemma}{Lemma}

\theoremstyle{thmstyletwo}
\theoremstyle{thmstylethree}

\input{MathSymbols}

\begin{document}
\title[Optimization on the Oblique Manifold for Sparse Simplex Constraints via Multiplicative Updates]
{Optimization on the Oblique Manifold for Sparse Simplex Constraints via Multiplicative Updates}

\author*[1]{\fnm{Flavia} \sur{Esposito}}\email{flavia.esposito@uniba.it}

\author[2]{\fnm{Andersen} \sur{Ang}}\email{andersen.ang@soton.ac.uk}

\affil*[1]{\orgdiv{Department of Mathematics}, \orgname{ University of Bari Aldo Moro}, \orgaddress{\city{Bari}, \postcode{70125}, \country{Italy}}}

\affil[2]{\orgdiv{School of Electronics and Computer Science}, \orgname{University of Southampton}, \orgaddress{\city{Southampton}, \postcode{SO17 1BJ}, \country{United Kingdom}}}

\abstract{Low-rank optimization problems with sparse simplex constraints involve variables that must satisfy nonnegativity, sparsity, and sum-to-1 conditions, making their optimization particularly challenging due to the interplay between low-rank structures and constraints.
These problems arise in various applications, including machine learning, signal processing, environmental fields, and computational biology.
In this work, we propose a novel manifold optimization approach to efficiently tackle these problems.
Our method leverages the geometry of oblique manifolds to reformulate the problem and introduces a new Riemannian optimization method based on Riemannian gradient descent that strictly maintains the simplex constraints.
By exploiting the underlying manifold structure, our approach improves optimization efficiency.
Experiments on synthetic and real datasets demonstrate the effectiveness of the proposed method compared to standard Euclidean and Riemannian methods, paving the way for broader applications.
}
 
\keywords{manifold, optimization, simplex, nonnegativity, sparsity}
%\pacs[JEL Classification]{D8, H51}
\pacs[MSC Classification]{15A23, 65K10, 49Q99, 90C26, 90C30}

\maketitle

\section{Introduction}\label{sec1}
Low-rank decomposition techniques, such as Nonnegative Matrix Factorization (NMF), exploit the fact that high-dimensional matrices can be well-approximated by a sum of rank-one components~\cite{gillis2020nonnegative}.
This reduces the computational complexity of the operations in the optimization while preserving the essential structure of the data~\cite{udell2019big}.
In NMF, a nonnegative matrix $\X\!\in\!\IRmn_+$ is approximated as the product of two low-rank matrices $\W\!\in\!\IRmr_+, \H\!\in\!\IRrn_+$, where $r \!\leq\! \min(m,n)$ is a pre-determined rank parameter.
The factor matrices from NMF can be obtained by solving the following minimization problem $\min_{\W\geq 0,\H\geq 0}\!\frac{1}{2}\|\X\!-\!\W\H\|_F^2$.
Additional constraints, such as sparsity or smoothing (to enforce temporal continuity or spatial coherence), are often imposed on the factors $\W$ or $\H$ to improve interpretability and performance \cite{del2025penalizing}.
{In many applications, the matrix $\W$ is assumed to be known or precomputed, for instance as a dictionary derived from prior information.
In this setting, the problem reduces to estimating the coefficient matrix $\H$, which typically encodes quantities such as abundances (in hyperspectral unmixing) or activations (in source separation and topic modeling) \cite{gillis2020nonnegative,ang2019algorithms,yang2010blind,qmatrix}.
This subproblem is therefore of independent interest and often involves additional structural constraints.}

\begin{figure}[h!]
\centering
\includegraphics[width=0.8\linewidth]{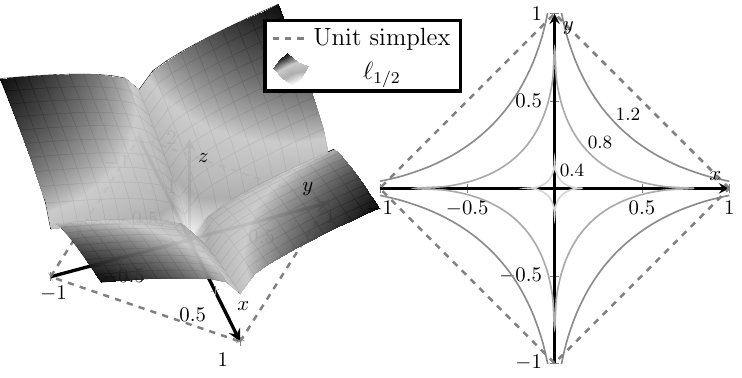}
\caption{The plot of the (nonconvex) $\ell_{1/2}-$quasi-norm and the simplex constraint.
{ Left: The simplex constraint is represented by the diagonal line intersecting the axes at $1$, visible as the dotted line in the left panel and the diagonal segment in the right panel.}
Right: the level set plot.
Minimizing the $\ell_{1/2}-$quasi-norm on the simplex gives the solution sitting on the corner of the simplex, leading to sparsity.
}
\label{fig:l0.5}
\end{figure}

In this work, we focus on the subproblem in $\H$ while keeping $\W$ fixed, 
proposing a novel formulation that differs from standard NMF approaches.
This choice is motivated by the fact that, in most applications, the interpretability of the factorization depends primarily on $\H$, where the relevant constraints are imposed.
We focus on combining three conditions: \textbf{nonnegativity}, \textbf{sparsity},  and \textbf{sum-to-1} into the subproblem in $\H$.
By enforcing these constraints on $\H$, the factorization acquires a probabilistic interpretation, enabling the evaluation of the relative contribution of each column of $\X$ (i.e., each sample or element) to the latent factors represented by $\H$.
These conditions on $\H$ finds applications:
\begin{itemize}
\item In biomedical contexts where $\X$ represents gene expression data (genes as rows and samples as columns), and $\H$ quantifies the relative contribution or abundance of latent biological factors such as metagenes across samples~\cite{brunet,Boccarelli2023}.
\item In remote sensing, such as hyperspectral unmixing, where columns of $\X$ contain spectral signatures of each pixel and $\H$ represents the fractional abundances of endmember materials; enforcing the constraints on $\H$ ensures physically meaningful interpretations of the contribution of each material within a pixel~\cite{ang2019algorithms,gillis2020nonnegative,settembre2025advancing}.
\end{itemize}
Specifically, we reformulate the optimization in $\H$ and solve it directly as a quartic optimization problem, where the nonnegativity and sum-to-1 constraints are intrinsically handled by a Hadamard parametrization on an oblique manifold.
As discussed later, this leads to an $\ell_1$-norm regularized quartic formulation.
To provide context, we first introduce the classical formulation of the subproblem in $\H$, which corresponds to the task of estimating the abundance or activation matrix $\H$ given a fixed, pre-computed dictionary matrix $\W$ and the observation matrix $\X$.
We consider the following Nonconvex-sparse simplex least squares \eqref{problem} problem
\begin{equation}\label{problem}
\argmin_{\H}
\tfrac{1}{2}\|\X-\W\H\|_F^2+
\lambda
\|\H\|_{1/2}
~~~~\st~
\underbrace{
\H \geq \zeros, 
~ \H^\top \ones_r = \ones_n,
}_{\text{columns of $\H$ in unit simplex}
}
\tag{NSSls}
\end{equation}
where $\|\cdot\|_F$ is the Frobenius norm,
$\lambda\!\in\!\IR_+$ is a penalty hyperparameter\footnote{{Hyperparameters are user-set parameters that affect models' behavior and performance, such as learning rate, architecture, regularization coefficients, and batch size, and must be tuned for optimal results \cite{esposito2023theoretical}; here, $\lambda$ controls the trade-off between the data fidelity term and the sparsity-inducing penalty.}}, 
$\ones_r\!\in\!\IR^r$ is vector-of-1,
and the nonconvex $\ell_{1/2}$-quasi-norm \cite{xu2012l}, for 
promoting sparsity\footnote{$\ell_{1/2}$ is $\ell_p$-norm with $p=1/2$.
This is a quasi-norm because it is not homogeneous.
}, 
is
\[
\textstyle
\|\H\|_{1/2}
=
\sum_{i=1}^r{\sum_{j=1}^n{\sqrt{H_{ij}}}}.
\tag{$\ell_{1/2}-$\text{quasi-norm}}
\]
{ This latter, $\|\H\|_{1/2}$, represents the $\ell_{1/2}$ quasi-norm, which evaluates the sum of the square roots of the absolute values of the entries.}
The $\ell_{1/2}$ quasi-norm can be seen as a non-convex approximation of the $\ell_0$ pseudo-norm that combines the sparsity-promoting effect of the element-wise $\ell_1$-norm with the smoothing properties of the $\ell_2$-norm.
By penalizing small coefficients more strongly and large ones more mildly, this quasi-norm drives weaker components toward zero more effectively, yielding sparser and more interpretable solutions.
This balance between the $\ell_1$ and $\ell_2$ behaviors makes the $\ell_{1/2}$ quasi-norm \cite{zeng2014l} particularly useful in applications such as compressive sensing and sparse regression, where both sparsity and stability are desired \cite{chartrand2007exact,chartrand2007nonconvex}.
{This quasi-norm with this specific power is used as it offsets the squaring effect of the Hadamard parametrization that will be introduced later.}
{As shown,} model \eqref{problem} encourages the columns in $\H$ to be concentrated on the boundary of a simplex  (defined by the nonnegativity and sum-to-1 constraints on each column), which promotes sparsity, while maintaining nonnegativity.
See Fig.~\ref{fig:l0.5}.

Problem \eqref{problem} has been addressed in \cite{qian2011hyperspectral} based on Euclidean geometry.
We reformulate the problem within a Riemannian framework, incorporating the constraints directly into the optimization process through fast and efficient  updates.

\textbf{Contribution.}
{In this work, we introduce a novel Riemannian approach for solving \eqref{problem}. The proposed approach consists of reformulating \eqref{problem} as a manifold optimization problem  on the oblique manifold and solving the resulting problem using an adaptation of Riemannian Multiplicative Update (RMU), an optimization method introduced in \cite{esposito2025new}. Thus, throughout the paper, RMU denotes the optimization algorithm, whereas the proposed approach refers to the combination of the manifold reformulation and the RMU solver. RMU is based on an approximate Riemannian gradient descent and employs a projection-free strategy that maintains the smoothness requirements on the manifold, namely the differentiability of the retraction and the continuity of the trajectory along the manifold geometry.}
The advantage of RMU is computational efficiency in handling the simplex constraint: RMU solves \eqref{problem} by incorporating the constraint directly into the minimization process. 
By reinterpreting the constrained problem in $\H$ as a manifold-based optimization, applying an efficient Riemannian solver, and providing robust experimental validation against Euclidean and alternative Riemannian methods, we demonstrate the significant practical benefits of this specific combination.
Moreover, the proposed approach differs from related existing works as follows:
\begin{itemize}
\item in \cite{qian2011hyperspectral}, we have a similar problem without the sum-to-1 constraint, as the objective function to be minimized is the same.
\item in \cite{ang2021fast} we consider a problem with a similar constraint but with an additional hyperplane constraint and a different objective function.
Specifically, the considered constraint is the k-capped simplex ($\Delta_k$), which can be interpreted as a hypercube intersected by a hyperplane.
\item in \cite{guo2021sparse}, we have the same problem, but it is solved with a different optimization approach based on the Riemannian conjugate gradient (RCG) algorithm.
\end{itemize}
We compare RMU with RCG and a variant of \cite{qian2011hyperspectral}, see Section~\ref{sec_exp} for detailed experimental comparisons.

\textbf{Paper Organization.}
Section\,\ref{sec_main} presents the manifold formalization of \eqref{problem}, reviews and describes manifold optimization for solving it.
{In this section we also provide a new local convergence theorem for the method proposed}.
Section\,\ref{sec_exp} presents numerical experiments.
Section\,\ref{sec_conc} concludes the paper.

\section{Manifold formulation}\label{sec_main}
Constrained problem \eqref{problem} can be viewed as an unconstrained optimization problem on a Riemannian manifold, where the constraints define a smooth search space without boundaries.
In this setting, the optimization is performed directly on the manifold, effectively eliminating the need for explicit constraint handling   \cite{boumal2023intromanifolds}.

The literature contains several approaches for handling constraints and penalizations similar to those in  \eqref{problem}~\cite{qian2011hyperspectral,guo2021sparse}.
However, many of these methods rely on Euclidean optimization and enforce the constraints through a post-normalization step, which may affect both efficiency and theoretical consistency.
In contrast, we formulate \eqref{problem} within the framework of Riemannian optimization on the oblique manifold, allowing the constraints to be intrinsically incorporated into the optimization process.
We also note that although the Riemannian Conjugate Gradient (RCG) method \cite{guo2021sparse} uses a similar formulation, it requires more steps and incurs higher per-iteration computational costs than the proposed method, as we will observe in Section \ref{sec_exp}.

\subsection{Formulation in oblique manifold}
Following \cite[Sec.III A]{guo2021sparse}, we reformulate the sum-to-1 constraints in \eqref{problem} by introducing a rectangular matrix $\A\in\IRrn$.
{We assume that $\H$ has the form} $\H\!=\!\A\!\odot\!\A$ {with $\odot$ the Hadamard product and}, we embed $\A$ in the rank-$r$ oblique manifold:
\[
\cO\cB(r,n) 
~=~
\{ \A\in\IRrn \,|\, \diag(\A^\top \A)=\I_n \},
\tag{Oblique manifold}
\]
where $\I_n$ is the identity, and $\diag:\IRnn\!\to\!\IRnn$ takes the diagonal elements of the input matrix to form a diagonal matrix\cite{absil2009optimization}.
We remark that the Hadamard reformulation $\H\!=\!\A\!\odot\!\A$, introduced in \cite[Section III A]{guo2021sparse}, is also recently considered in \cite{li2023simplex} with more detailed analysis.
We present the following lemma, {inspired by \cite[Section III A]{guo2021sparse}.
While \cite{guo2021sparse} introduced the Hadamard parametrization, we provide the explicit formal proof of its bidirectional equivalence (as stated in Lemma~\ref{lem:HA}) for completeness.}
\begin{lemma}\label{lem:HA}
Columns of $\H = \A \odot \A$ stay in the simplex if and only if $\A\in\cO\cB(r,n)$.
\begin{proof}
Nonnegativity is given by the definition $\H\!=\!\A\!\odot\!\A\!\geq\!\zeros$.
The sum-to-1 constraint naturally follows: for all $j$ we have $\sum_{i=1}^r \! H_{ij}
\!=\!
\sum_{i=1}^r \! A_{ij}^2
\!=\!
(\A^\top\A)_{jj}$.
Then by $\A\!\in\!\cO\cB(r,n)$, we have this sum equals to $1$.
\end{proof}
\end{lemma}
\begin{remark}[\textbf{Uniqueness.}]
The factor matrix $\A$ in  Lemma~\ref{lem:HA} is non-unique.
For example we can have $\B = -\A$ giving $\H = \B\odot\B = (-\A)\odot(-\A) = \A \odot \A$.
\end{remark}
Based on Lemma~\ref{lem:HA}, {we reformulate \eqref{problem} as a manifold optimization problem, introducing a new function $f:\IRrn\!\to\!\IRrn$, which we refer to as the L1-quartic least squares problem \eqref{riemanianfunc_4}}
\begin{equation}\label{riemanianfunc_4}
\argmin_{\A\in\cO\cB(r,n)}
\Big\{
f(\A)
=
\tfrac{1}{4}
\|\X-\W(\A\odot\A)\|_F^2
+\lambda\|\A\|_1
\Big\},
\tag{L1-quartic-ls}
\end{equation}
where $\|\A\|_1\!=\!\sum_{ij}|A_{ij}|$ is the entry-wise $\ell_1$-norm of $\A$ due to $(x^\frac{1}{2})^2\!=\!|x|$.

\begin{remark}[\textbf{Why $\ell_{1/2}$-quasi-norm.}]
Other nonconvex sparsity-inducing norms such as the $\ell_1$--$\ell_2$~\cite{yin2015minimization} are also available.
We use the nonconvex $\ell_{1/2}$-quasi-norm due to Lemma~\ref{lem:HA}: by $\H\!=\!\A\!\odot\!\A$, the Hadamard product cancels with the $1/2$-power in the $\ell_{1/2}$-quasi-norm.
\end{remark}

\begin{remark}[Nonconvexity in objective]
Both the objective functions in \eqref{problem} and \eqref{riemanianfunc_4} are nonconvex.
In \eqref{problem}, the term $\ell_{1/2}\text{-quasi-norm}$ is nonconvex.
In \eqref{riemanianfunc_4}, the term $\|\X-\W(\A\odot\A)\|_F^2$ is quartic in $\A$ and generally nonconvex \cite{ahmadi2013np}.
{ This highlights a challenge: replacing the $\ell_{1/2}$-quasi-norm with the $\ell_1$-norm via the Hadamard parametrization ($\H = \A \odot \A$) does not make the overall problem convex.
The quartic data fidelity term ensures that the cost landscape remains highly non-convex, meaning optimization algorithms may be susceptible to local minima.}
\end{remark}

Since the solution to the manifold optimization problem \eqref{riemanianfunc_4} also solves \eqref{problem}, (via Lemma~\ref{lem:HA}), we briefly review key concepts from Riemannian optimization to solve \eqref{riemanianfunc_4}.
We keep the material minimum, for details see~\cite{boumal2023intromanifolds, absil2009optimization,
guo2021sparse}.

\subsection{Oblique manifold}
The oblique manifold $\cO\cB(r,n)$ is the set of matrices {whose} columns have unit norm.
The set $\cO\cB(r,n)$ is also the product of $r$ spheres
$\cS^{n-1}\!=\!\{\x\!\in\!\IRn |\, \|\x\|_2 \!=\!1\}$, that is, from a geometric point of view, each column of those matrices in $\cO\cB(r,n)$ lies on individual spheres.
We refer the reader to \cite{Tu2010em,boumal2023intromanifolds,absil2009optimization} for details on the oblique manifold.
Our optimization problem can be seen as {unconstrained optimization problem} over a curved space, and we make use of Riemannian optimization tools to solve \eqref{riemanianfunc_4}, where each update direction must preserve the column norms, maintaining a smooth but non-linear geometric structure.
To this end, several fundamental concepts from Riemannian optimization are required: 
\begin{itemize}
    \item the tangent space, which provides a local linear approximation of the manifold;
    \item  the projection onto the tangent (and normal) spaces, used to compute the Riemannian gradient that satisfies the geometric constraints, unlike the Euclidean gradient, which is generally not tangent to the manifold;
    \item and the retraction, which maps an updated point back onto the manifold after each optimization step.
\end{itemize}
Below, we detail these concepts for the oblique manifold.
The tangent space of $\cO\cB(r,n)$ at a point $\A$, denoted as $\cT_{\A} \cO\cB(r,n)$, {is defined as}
\[
\cT_{\A} \cO\cB(r,n) 
= 
\{
\Z\in\mathbb {R}^{r\times n} \,|\, \diag(\A^\top \Z) = \zeros
\},
\]
with the associated null space $\cN_{\A}\cO\cB(r,n)\!=\!\{ \A\D \,|\, \D\!\in\!\IRnn \text{ is diagonal}\}$.
The projectors over these spaces at $\A$ for a generic matrix $\Z\in\IRrn$ are
\begin{equation}\label{projectors}
\cP_{\cT_{\A}}(\Z) = \Z-\A \diag(\A^\top\Z),
\quad
\cP_{\cN_{\A}}(\Z) = \A \diag(\A^\top\Z).
\end{equation}
Lastly, we use metric retraction to move a point back to the manifold { from the ambient Euclidean space}:
\begin{equation}\label{metricRetraction}
\cR_{\A}(\Z) 
= 
(\A+\Z)\big(
\diag \big(
(\A+\Z)^\top(\A+\Z)
\big)^{-1/2}
\big).
\end{equation}
{Before deriving the new RMU for
\eqref{riemanianfunc_4}, we recall the RMU method proposed in \cite{esposito2025new} and its formulation, as well as the new convergence results for the oblique manifold}.

\subsection{Riemannian Multiplicative Update on the Oblique Manifold}\label{sec:2:subsec:RMU}
In this section we specialize RMU to the oblique manifold $\cO\cB(r,n)$. Throughout this section we assume that optimization is performed over $\cO\cB(r,n)$.
RMU is a method proposed in \cite{esposito2025new} based on Riemannian gradient descent~\cite{boumal2023intromanifolds}, to solve a nonnegative problem $\argmin_{\x\in\cM} f(\x)$, with $\x\!\geq\!0$.
{ In our context referring to \eqref{riemanianfunc_4} problem, the generic variable $\x$ corresponds to the factor matrix $\A$, the manifold $\cM$ is the oblique manifold $\cO\cB(r,n)$, and the objective function $f$ is the L1-quartic least squares function cost.}
In our setting, RMU employs the sign-wise splitting
\begin{equation}\label{sign_split}
\grad f(\A)\!=\!\grad^+ f(\A)-\grad^- f(\A),
\end{equation}
to be explained later, 
ii) a metric retraction $\cR(\cdot)$ {defined in \eqref{metricRetraction}}, 
iii) and an appropriate stepsize  $\alpha$ {that we detail in the following}.
{
According to this, we can state the following corresponding theorem for RMU in the case of $\cO\cB(r,n)$\cite{esposito2025new}.
\begin{theorem}
\label{prop:RMU_update}
Let $\IRrn$ be the ambient Euclidean space of a manifold $\cO\cB(r,n)$ and consider $f:\cO\cB(r,n)\!\to\!\IR$ be a continuously differentiable function bounded from below.
Let $\V^k\!=\!-\grad f(\A^k)$ be the anti-parallel direction given by the Riemannian gradient of $\cO\cB(r,n)$ at $A^k$, and let $\cR_{A^k}$ be the metric retraction onto $\cO\cB(r,n)$.
If a nonnegative $\A^k$ is updated by $\A^{k+1} = \cR_{\A^k}( \balpha^k\!\odot\!  \V^k)$
with an element-wise\footnote{{The symbol $\oslash$ is used to indicate the element-wise division.}} stepsize $\balpha^k\!=\!\A^k\!\oslash\!\grad^+ f(\A^k)\!\in\!\IRrn$with the sign-wise splitting $\grad^+ f$ in \eqref{sign_split}, then $\A^{k+1}\!\geq\!\zeros$ and it is still on $\cO\cB(r,n)$.
\end{theorem}}

\paragraph{Element-wise operations may not preserve the tangent condition}
We note that RMU is not strictly a Riemannian gradient descent method, as the element-wise operations, as linear transformations in general, do not preserve the tangent condition of a manifold at a point.
However, we can still perform metric retraction regardless of the point is on the tangent space or not.
This phenomenon may contribute to a slower convergence rate, a tendency that was clearly reflected in the outcomes of our experimental analysis in Section \ref{sec_exp}.
{ Furthermore, without strict adherence to the tangent space, the rigorous theoretical guarantees of monotonic descent standard in Riemannian gradient methods are not assured, meaning the algorithm relies heavily on the local metric preservation properties of the retraction to achieve convergence in practice.
In the case of the oblique manifold, the convergence of RMU is not affected by the lack of a tangent condition, since the update is mapped back to the manifold through the retraction, and a descent property can still be established locally.
To this aim, we have the following local convergence theorem.
%\begin{theorem}[Local convergence]\label{thm:localconvergence}
%Let $\A^{k+1} = \cR_{\A^k}(-\boldsymbol{\alpha}^k\odot \grad f(\A^k))$ be the RMU update as in Theorem \ref{prop:RMU_update},  where each component of $\balpha^k$ is bounded in $[\alpha_{min}, \alpha_{max}]$ with $\alpha_{min}>0$.
%For sufficiently small $\alpha_{max}$, the sequence $\{\A^k\}_k$ generated by RMU satisfies $\lim_{k \rightarrow \infty}\|\grad f(\A^k)\|= 0$.
%\end{theorem}

\begin{theorem}[Local Convergence]\label{thm:localconvergence}
Let $f:\cO\cB(r,n)\!\to\!\mathbb{R}$
be a continuously differentiable function bounded from below in an open neighborhood containing the sequence $\{\A^k\}_{k\ge0}$ generated by RMU. Assume that: i) the Riemannian gradient $\grad f$ is locally Lipschitz continuous; ii) the retraction is the metric retraction $\cR$ defined in~\eqref{metricRetraction}; iii) each component of the stepsize matrix satisfies
\[
0<\alpha_{\min}\le (\balpha^k)_{ij}\le \alpha_{\max},
\]
where $\alpha_{\max}$ is sufficiently small;
iv) either the iterates avoid the nondifferentiable points of the entrywise $\ell_1$-norm, or $f$ is replaced by a locally smoothed approximation.
Let the iterates be generated by RMU as $\A^{k+1}
=
\cR_{\A^k}
\big(
\!-\!\balpha^k \!\odot\! \grad f(\A^k)
\big).
$
Then there exists a constant $c>0$ such that
$f(\A^{k+1})
\!\leq\!
f(\A^k)
-
c\,
\|\grad f(\A^k)\|^2,
$
for all sufficiently large $k$. Consequently, $
\sum_{k=0}^{\infty}
\|\grad f(\A^k)\|^2 \!<\! \infty
$,
and
$
\lim_{k\to\infty}
\|\grad f(\A^k)\| \!=\! 0
$.
Hence every accumulation point of the sequence $\{\A^k\}$ is a stationary point of problem~\eqref{riemanianfunc_4}.
\end{theorem}
\begin{proof}
We prove that the descent property holds for RMU by working with a smooth extension of the retraction $\cR_{\A}$ to a neighborhood of the ambient space.
In the following, we work on the first-order expansion of the retraction showing that only the tangential component affects the update, ensuring that the descent analysis remains valid.
\\
Let $\A \in \cO\cB(r,n)$, which implies $\|\A_j\|\!=\!1$ for all columns $j$.
Consider the retraction for a \textbf{small} (possibly non-tangent) $\Z$: we have $\cR_{\A}(\Z) =
(\A+\Z)\big(\diag\big((\A+\Z)^\top(\A+\Z)\big)\big)^{-1/2}$ from $\eqref{metricRetraction}$.
As $\|\A_j\|\!=\!1$, we expand $\|\A_j\!+\!\Z_j\|\!=\!\sqrt{1\!+\!2\A_j^\top \Z_j\!+\!\|\Z_j\|^2}$.
By the Taylor approximation $(1+\epsilon)^{-1/2} = 1-\tfrac{1}{2}\epsilon+O(\epsilon^2)$, we get $\tfrac{1}{\|\A_j+\Z_j\|}
=
1 - \A_j^\top \Z_j + O(\|\Z_j\|^2)$.
Substituting this back into the column-wise retraction gives:
\[
\begin{array}{rcll}
\begin{bmatrix}
\cR_{\A}(\Z)   
\end{bmatrix}_j
&=& \A_j+\Z_j/\|\A_j+\Z_j\|
\\
&=& (\A_j+\Z_j)\big(1 - \A_j^\top \Z_j + O(\|\Z_j\|^2)\big)
\\[1mm]
&=& \A_j + \Z_j - \A_j(\A_j^\top \Z_j) + O(\|\Z_j\|^2).
\end{array} 
\]
By the definition of the tangent space and the corresponding projection $\cP_{\cT_{\A}}(\Z) = \Z-\A \diag(\A^\top\Z)$ given in \eqref{projectors}, we see that $\Z_j - \A_j(\A_j^\top \Z_j) = \cP_{\cT_{\A_j}}(\Z_j)$.

Thus, the first-order Taylor expansion of the retraction in matrix form is simply: 
\[
\cR_{\A}(\Z) =\A + \cP_{\cT_{\A}}(\Z) + O(\|\Z\|^2).
\]  
Now, consider the update $\A^{k+1}=\cR_{\A^k}(\boldsymbol{\eta}^k)$  with $\boldsymbol{\eta}^k=-\boldsymbol{\alpha}^k\odot \grad f(\A^k)$.
By assumption (iv), f is locally smooth.
The Taylor expansion of the cost function $f(\cR_{\A^k}(\boldsymbol{\eta}^k))= f\big(
\A^k + \cP_{\cT_{\A^k}}(\boldsymbol{\eta}^k) + O(\|\boldsymbol{\eta}^k\|^2) \big)$ gives:
\[
\begin{array}{rcll}
f(\cR_{\A^k}(\boldsymbol{\eta}^k))
&=& 
f(\A^k) + \langle \nabla f(\A^k), ~ \cP_{\cT_{\A^k}}(\boldsymbol{\eta}^k) \rangle 
+ O(\|\boldsymbol{\eta}^k\|^2)
\\[1mm]
&=&
f(\A^k) + \langle  \cP_{\cT_{\A^k}}(\nabla f(\A^k)), ~ \boldsymbol{\eta}^k \rangle 
+ O(\|\boldsymbol{\eta}^k\|^2)
\\[1mm]
&=&
f(\A^k) + \langle \grad f(\A^k), ~ \boldsymbol{\eta}^k \rangle 
+ O(\|\boldsymbol{\eta}^k\|^2), 
\end{array}
\]
where the last two equalities hold by the properties of the orthogonal projector and the definition of the Riemannian gradient.
Now, since $\langle \grad f(\A^k),\boldsymbol{\eta}^k\rangle \leq -\alpha_{\min}\|\grad f(\A^k)\|^2$ and $\|\boldsymbol{\eta}^k\|^2 \leq \alpha_{\max}^2\|\grad f(\A^k)\|^2$, we have:
\[
f(\A^{k+1})
\leq
f(\A^k)
-
c\|\grad f(\A^k)\|^2,
\qquad c = \alpha_{\min}-S\alpha_{\max}^2
\] 
where $S$ is the existing nonnegative constant such that $\|O(\|\eta\|^2)\|\!\leq\! S\|\eta\|^2$ .
For $\alpha_{\max}$ sufficiently small we have $c>0$, which guarantees a descent direction.
As $f$ is bounded from below, the sequence $\{f(\A^k)\}_k$ is decreasing and convergent, which implies $\|\grad f(\A^k)\| \!\to\! 0$.
\end{proof}
}
{We are now ready to describe the details of RMU for solving \eqref{riemanianfunc_4}.}

\subsection{Algorithm for minimizing  \eqref{riemanianfunc_4} over $\cO\cB(r,n)$}
Under the assumptions of Theorem 1, the objective function in \eqref{riemanianfunc_4} satisfies the hypotheses of RMU. 
{In fact, although the entry-wise $\ell_1$-norm is not differentiable on the
entire ambient space, the RMU iterates remain nonnegative when initialized
with $A^0 \geq 0$, as established in Theorem~1. Therefore, along the RMU iterates, the $\ell_1$-norm reduces to the linear function $\|\A\|_1=\sum_{i,j}{A_{ij}}$, which is continuously differentiable on the nonnegative region.
Thus, Theorem~1 can be applied to this smooth representation of the objective
on the invariant nonnegative region.}
Therefore, the update
$
\A^{k+1}
\!=\!
\cR_{\A^k}\big(\!
-\!\balpha^k\!\odot\! \grad f(\A^k)
\big)
$
converges to a stationary point of function \eqref{riemanianfunc_4} for $\A\in \cO\cB(r,n)$.
In this update, function $\cR(\cdot)$ is the retraction \eqref{metricRetraction}, and $\grad f(\cdot)$ is the Riemannian gradient computed with the orthogonal projector over the tangent space in \eqref{projectors} over the Euclidean (sub-)gradient:
\[
\nabla f(\A)
~=~ 
-
2\W^\top\big(\X-\W(\A\odot\A)\big)\odot \A
+
\lambda' \sgn(\A),
\]
where $\sgn$ is the element-wise sign function (subgradient of the $\ell_1$-norm) with $\sgn(x) = 1$ if $x>0$ and $\sgn(x)=-1$ if $x<0$ and $\sgn(x)\in [-1,1]$ if $x=0$, see for details~\cite[Example 3.4]{beck2017first}.
The symbol $\lambda'$ is a rescaled $\lambda$ (with respect to the $1/4$ in $f$ when taking the derivative).

Let $\Q=\W^\top\W$ and $\P=\W^\top\X$, the Riemannian (sub-)gradient is defined as:
\begin{equation}\label{gradient}
\hspace{-2mm}
\begin{array}{lll}
\grad f(\A) 
=
\underbrace{
\big(\Q(\A\odot\A)\big)\odot\A
+
\A \diag\big[\A^\top (\P\odot\A )\big]
+ \lambda'
\sgn(\A)
}_{\grad^+ f} 
\\
\qquad\qquad~
- 
\underbrace{
\biggl\{
\P\odot\A
+
\A \diag\big[\A^\top\big(\Q(\A\odot\A)\big)\odot\A\big]
+
\lambda'
\A \diag\big[\A^\top \sgn(\A)\big]
\biggr\}
}_{\grad^- f},
\end{array}
\end{equation}
where we used a sign-wise splitting 
$\grad f(\A)\!=\!\grad^+ f(\A)-\grad^- f(\A)$ as in \eqref{sign_split}.
Then, from Theorem~\ref{prop:RMU_update}, considering $\V^k\!=\!-\grad f(\A^k)$, the sign-wise splitting in \eqref{gradient} and the element-wise stepsize $\balpha^k$, we obtain the RMU update as:
\begin{equation}\label{eqn:updt_BA}
\B^k 
\!=\!
\A^k\odot \grad^- f(\A^k)\oslash \grad^+f(\A^k),    
\quad 
\A^{k+1}
\!=\!
\B^k
\oslash 
\diag\big[(\B^k)^\top \B^k \big]^{1/2},
\end{equation}
where $\oslash$ is element-wise division and, $\B^k$ represents an intermediate Euclidean step taking the element-wise scaled gradient direction before retraction, and $\balpha^k$ is implicitly embedded in the ratio $\A^k \oslash \grad^+ f(\A^k)$.
Finally,  after $K$ iterations, we have $\A^K$, and we compute $\H\!=\!\A^K\!\odot\!\A^K$.

Algorithm~\ref{alg:RMU} summarizes in pseudo-code the previous steps for solving \eqref{riemanianfunc_4} with RMU on the oblique manifold.
\begin{algorithm}[h!]
\label{algo_ob}
\caption{{ RMU solver for \eqref{riemanianfunc_4}}}
\label{alg:RMU}
\begin{algorithmic}[1]
\Require Data $\X\in\IRmn_+$, rank $r$, $\W\in\IRmr_+$, $\lambda > 0$
\State Initialize $\A^{0} \in \mathcal{OB}(r,n)$
\For{$k = 1$ to $K$}
\State Compute $\grad^+ f(\A^k)$ and $\grad^- f(\A^k)$ as in~\eqref{gradient} considering $\A = \A^k$.
\State Update $\B^k, \A^{k+1}$ as in~\eqref{eqn:updt_BA}
\EndFor
\State \Return $\H = \A^{K} \odot \A^{K}$ 
\end{algorithmic}
\end{algorithm}

\begin{remark}{
The entry-wise $\ell_1$-norm is generally nonsmooth at matrices containing
zero entries. However, RMU is initialized with a nonnegative matrix and
preserves nonnegativity at every iteration. Hence, along the RMU iterates,
\[
    \|A\|_1 = \sum_{i,j} A_{ij},
\]
which is the restriction to the nonnegative orthant of a continuously
differentiable linear function. Its Euclidean gradient is the all-ones
matrix. At an entry such that $A_{ij}=0$, this corresponds to selecting
\[
    1 \in [-1,1] = \partial |A_{ij}|,
\]
which is a valid element of the subdifferential of the original
$\ell_1$-norm. Therefore, the gradient employed by RMU is consistent both
with the smooth representation of the objective on the nonnegative region and with the subdifferential of the original objective.}
The subgradient of $\|\cdot\|_1$, denoted as $\partial \|\A\|_1$, has the form
\[
\partial \|\A\|_1 
= \sgn(\A) = \G \in \IRrn,
~~ \text{where }
G_{ij} \in \begin{cases}
   \{0\} & \text{ if }A_{ij} \neq 0,
   \\
   [-1,1] & \text{ if }A_{ij} = 0.
\end{cases}
\]
Moreover, we remark that the structure $\sgn(\A)-\A \diag\big[\A^\top\sgn(\A)\big]$ comes from the projectors in \eqref{projectors} applied on $\partial \|\A\|_1 $.
\end{remark}

\begin{remark}
$\A^k\!\oslash\! \grad^+f(\A^k)$ in \eqref{eqn:updt_BA} is implemented as 
$\A^k\!\oslash\! \big(\grad^+f(\A^k) + \epsilon\big)$, with $\epsilon\!>\!0$ a small positive constant preventing division by zero introduced to avoid division by zero.
\end{remark}

\paragraph{Validity of the Retraction of Point Outside the Tangent Space}
As mentioned in Section\,\ref{sec:2:subsec:RMU} the element-wise operations may not preserve the tangent condition, which is generally required for the retraction $\cR(\cdot) : \cT_{\A} \cO\cB(r,n) \to \cO\cB(n,r)$.
We emphasize that because the retraction $\cR$ is specifically designed to normalize the columns of the input matrix, it maps points back to $\cO\cB(r,n)$, regardless of whether the input matrix lies in the tangent space.
In fact, the normalization operator defining the retraction can also be applied to points
$\Z \notin \cT_{\A} \cO\cB(r,n)$.
Furthermore, depending on the location of the point $\Z$, the metric properties of the retraction $\cR(\Z)$ are approximately preserved for small $\Z$, which is sufficient in practical optimization.
This implies that the retraction remains well-defined when applied to directions outside the tangent space and, by its first-order consistency, we can assert that locally the descent property and convergence are guaranteed as stated in Theorem~\ref{thm:localconvergence}.

\section{Numerical Experiments}\label{sec_exp}
In this section, we show numerical results of RMU compared with other methods and tested on different synthetic and real datasets.
Specifically, we first evaluate the proposed approach on synthetic datasets against Euclidean and Riemannian methods.
Second, we present a preliminary experiment on a real dataset to compare the proposed method strictly within the Riemannian optimization framework\footnote{The experiments are conducted in MATLAB 2023a on a machine with OS Windows 11 Pro on an AMD Ryzen 7 7700 8-Core Processor CPU 3.80 GHz, and 64GB RAM.}.

\textbf{Comparative methods.}
We compare RMU with three methods: one Riemannian optimization method and two standard Euclidean optimization methods.

\begin{itemize}
\item \textbf{Riemannian method}:
\cite{guo2021sparse} proposed a Riemannian Conjugate Method (RCG) for solving \eqref{problem}.
In addition to the Riemannian gradient, RCG also computes: i) a conjugate direction, a hybrid (Dai-Yuan and Hestenes-Stiefel) stepsize using vector transport, iii) update stepsize using Wolfe line search (with at most $5$ iterations, following \cite{gillis2020nonnegative}).

\item \textbf{Euclidean-based methods.}
We consider two Euclidean-based approaches.
\begin{itemize}
\item \textbf{EMUproj heuristic}:
For Problem~\eqref{problem} without the sum-to-1 constraint, the following update formula for $\H$ is proposed by 
\cite{qian2011hyperspectral}
\[\textstyle
\H_{k+1}^{\text{EMU}} 
= 
\H_k \odot \W^\top \M \oslash \Big(
\W^\top \W \H_k \oplus \tfrac{\lambda}{2} \sum_{ij}\sqrt{H_{ij}}
\Big),
\tag{EMU}
\]
which is a modified version of the work by Hoyer \cite{hoyer2002non}, where $\oplus$ is element-wise addition.
To solve Problem~\eqref{problem} by EMU, we consider the  heuristic
\[
\H_{k+1} = \proj_{\Delta}[ \H_{k+1}^{\text{EMU}}],
\tag{EMUproj}
\]
where $\proj_{\Delta}$ is the column-wise projection onto the unit-simplex.
Given a matrix $\X=[\x_1,\ldots,\x_n]$ with $n$ columns, then $\proj_{\Delta}(\X)$ replaces each $\x_i$ by $\x_i/\|\x_i\|_1$.

\item \textbf{SparseMU-proj heuristic}:
The sparse multiplicative update proposed by \cite{eggert2004sparse} has a form similar to EMU 
\[
\H_{k+1}^{\text{SMUL1}}= 
\H_k 
\odot 
\W^\top \M \oslash \Big(
\W^\top \W \H_k \oplus \lambda
\Big),
\tag{SMUL1}
\]
which minimizes a different objective function from the one in~\eqref{problem}, as
\[\textstyle
\argmin_{\H\geq 0}
\tfrac{1}{2}\|\X-\W\H\|_F^2+
\lambda
\|\H\|_{1}.
\]
We note that the difference between EMU and SMUL1 is on the regularization parameter, where the latter uses a constant value of $\lambda$, while the former adopts a dynamically changing $\lambda$ scaled by $\sum \sqrt{H_{ij}}$.
To solve Problem~\eqref{problem} by SMUL1, we consider the  heuristic with the same projection of EMUproj:
\[
\H_{k+1} = \proj_{\Delta}[ \H_{k+1}^{\text{SMUL1}}].
\tag{SparseMU-proj}
\]
We emphasize that SparseMU-proj is not designed to solve \eqref{problem}, however we still plot in our experiment for comparisons due to its proximity to our problem.

\begin{remark}
{ Note that applying these heuristic projections a posteriori destroys the monotonic decrease guarantees of the original Euclidean multiplicative update rules.
Consequently, these modified Euclidean baselines lack formal convergence guarantees, serving purely as heuristic benchmarks.}
\end{remark}
\end{itemize}
\end{itemize}
Table~\ref{tab:complexity} shows the per-iteration computational cost of the methods.
RMU has a lower per-iteration cost compared with RCG.
RMU {has a comparable cost with the modified Euclidean methods.}

\begin{itemize}
\item RCG solves NSSls using conjugate gradient, which involves additional computations (e.g., vector transport and line search).
In contrast, the Euclidean methods EMU and SparseMU, like our approach, use multiplicative updates and thus present comparable complexity.

\item The two Euclidean methods were originally proposed to enforce only sparsity.
In this work, for a fair comparison, we added a projection onto the unit simplex to include the sum-to-1 constraint.
\end{itemize}

\begin{table}[h!]
\centering
\caption{Per-Iteration complexities}
\label{tab:complexity}
\begin{tabular}{l|cccc}
\textbf{Operation} & \textbf{EMUproj} & \textbf{SMUL1} & \textbf{RMU} & \textbf{RCG} 
\\
\midrule
Element-wise $\sqrt{H}$ & $rn$ & -- & $rn$ & $rn$ 
\\
Column sums $\sum\sqrt{H}$ & $rn\!-\!1$ & -- & -- & -- 
\\
Sign computation & -- & -- & $rn$ & --
\\
Matrix mult. $W^\top WH$ & $(2r^2\!-\!r)n$ & $(2r^2\!-\!r)n$  & $(2r^2\!-\!r)n$ & $(2r^2\!-\!r)n$ 
\\
Diagonal operations & -- & -- & $6rn$ & -- 
\\
Element-wise operations & $4rn$ & $4rn$ & $10rn$ & $8rn$  
\\
Riemannian gradient & -- & -- & $4rn$ & $(4r^2\!+\!6r)n$ 
\\
Retraction & -- & -- & $3rn$ & $5rn$  
\\
Vector transport & -- & -- & -- & $4rn$  
\\
Wolfe line search & -- & -- & -- & $(10r^2\!+\!50r)n$  
\\
Projection (normalization) & $(3r\!-\!1)n$ & $(3r\!-\!1)n$ & -- & -- 
\\
\midrule
Total FLOPs & $(2r^2\!+\!8r\!-\!1)n$ & $(2r^2\!+\!6r\!-\!1)n$ & $(2r^2\!+\!24r)n$ & $(16r^2\!+\!74r) n$
\end{tabular}
\end{table}

\textbf{Experimental settings.}
We generate a matrix $\X$ as the product of two low-rank matrices $\W_{\text{true}},\H_{\text{true}}$ such as $\X\!=\!\W_{\text{true}}\H_{\text{true}}\!+\!\sigma \E$.
Both matrices $\W_{\text{true}}, \H_{\text{true}}$ are generated by sampling a uniform distribution, with $\H_{\text{true}}$ explicitly constrained to have exactly $s\%$ non-zero entries.
The matrix $\E$ is a noise matrix following a uniform distribution, with $\sigma\!\in\!\IR$ as scaling factor.
We perform 100 Monte Carlo runs of the algorithms according to the randomness of the dataset generation.
{ Each of the Monte Carlo runs uses a different random seed to regenerate the synthetic matrices ($\W_{\text{true}}$ and $\H_{\text{true}}$)  from the specified uniform distributions and different levels of noise according to the scaling factor ($\sigma=0,0.1,0.3$).
This has been done to ensure statistical robustness of the experiments.}

All the algorithms start with the same initialization and penalty parameter $\lambda$.
We initialize the algorithms by NNSVD \cite{boutsidis2008svd}.
Given an initial variable $H_{\text{ini}}$, we select the parameter $\lambda$ { scaled by} the ratio
$\|\X-\W\H_{\text{ini}}\|_F^2 / \|\H_{\text{ini}}\|_{1/2}$, which balances the least squares part and the sparsity part at the beginning of the optimization process.
We stop the algorithms when they reach either the maximum number of iterations $k_{\max}$ or the maximum runtime $T_{\max}$.
{ For the experiments, we strictly relied on these maximum iterations and runtime limits as the stopping criteria, rather than a relative tolerance on the objective, to evaluate the performance profile over the entire allocated time budget.}

\begin{figure}[h!]
\centering
\includegraphics[width=0.9\textwidth]{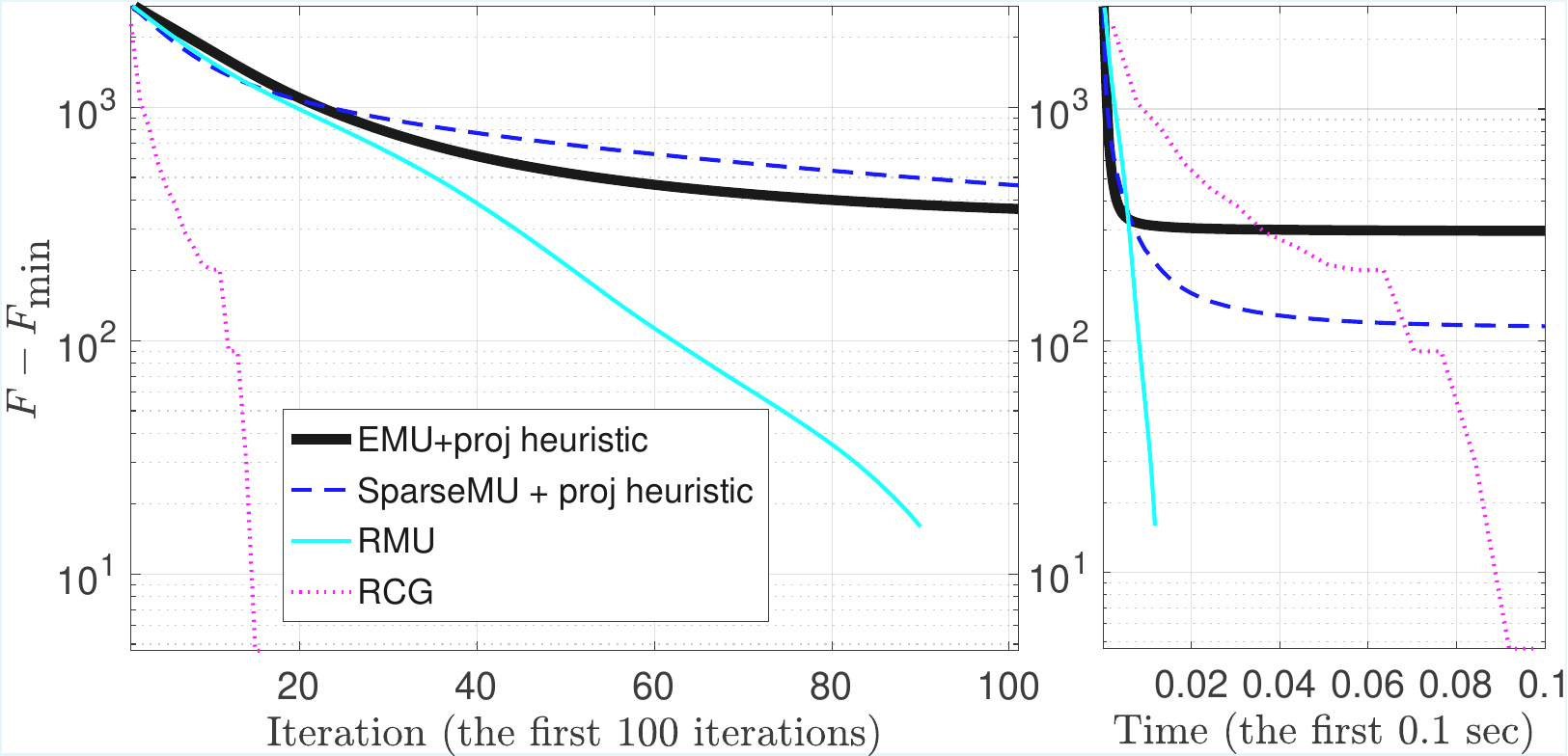}
\caption{A convergence plot of the methods { (where $F_{\text{min}}$ is the minimum objective value achieved across all compared algorithms for that specific run)}.
Here, in terms of iterations, the convergence of RCG has the best performance; however, it takes the highest amount of computational time for performing one iteration.
}
\label{fig:results}
\end{figure}

\textbf{AUC evaluation metrics.}
Because the methods have different per-iteration complexities, plotting the function value $F$ against the iteration $k$ or time $t$ alone provides an incomplete and unfair comparison.
See Fig.~\ref{fig:results}.

For each method we report the mean and standard deviation of the objective function, a ranking based on the area under the curve (AUC) metric, and a median sparsity measure of $\H$.
The ranking based on the AUC metric is computed over 100 Monte-Carlo runs.
The $i$th entry of the ranking vector records how many times an algorithm produced the $i$th best solution. 
For instance, a ranking tuple $(a, b, c, d)$ means the method ranked 1st $a$ times, 2nd $b$ times, 3rd $c$ times, and 4th $d$ times out of the 100 runs.
For each run, we execute all algorithms using the same initialization and modeling parameters, and obtain four sequences of the form $\{t_k, F_k\}_{k \in \IN}$, where $F_k$ is the objective value at iteration $k$, and $t_k$ is the corresponding runtime.
Because different algorithms have different per-iteration costs, we use AUC to compare their performance fairly. 
This time-based AUC evaluation metric is a customized metric proposed in this work to explicitly account for the varying per-iteration computational costs of different algorithms. 
Over the time interval $[0, T_{\max}]$, we define the AUC as $\int_0^{T_{\max}}\!F(t) \text{d}t$,
which represents the total accumulated objective value over time.
For a minimization problem, a smaller AUC indicates that the algorithm kept the objective low for most of the runtime.
Numerically, we treat $F(t)$ as a continuous function obtained by interpolating the points $\{t_k, F_k\}_{k \in \IN}$, and approximate the integral using the trapezoidal rule
$\Delta_{k}\!=\!(F_{k-1}+F_{k})(t_k-t_{k-1})/2$.
The AUC of a method is the sum of all $\Delta_k$.
We use AUC for the following reasons:
\begin{itemize}
\item It ignores the iteration index $k$ and evaluates performance in real time $t$, which allows a fair comparison of algorithms with different iteration counts.

\item The per-iteration computational cost is naturally captured in $t_k$; expensive iterations are beneficial only if they achieve a sufficiently large decrease in $F$.

\item When algorithms generate objective values at different timestamps $t_k$, interpolation creates comparable curves $F(t)$ over the common interval $[0, T_{\max}]$.
Thus, AUC enables fair comparison without requiring explicit synchronization of the curves.

\item AUC is robust to early stopping: if an algorithm terminates before $T_{\max}$, interpolation/extrapolation allows evaluation over the entire time window, ensuring consistent comparison across methods.
 
\end{itemize}

\textbf{Sparsity measure}: We also report the percentage sparsity of the output $\H_{\text{Last}}$ as $\text{nnz}(\H_{\text{Last}}\!>\!10^{-9})/nr \%$.
That is, we take entries in $\H_{\text{Last}}$ below $10^{-9}$ as zero, and count the zeros in $\H_{\text{Last}}$ with respect to the total number of entries $rn$.

\textbf{Results and Discussion.}
Now we report the results on three different-sized datasets, with the same sparsity level for $\H$ (set as $0.6$), and $3$ different noise levels $\sigma \in \{0,0.1,0.3\}$.
The datasets we refer to are considered as: i) small size data, $(m,n,r)\!=\!(20,100,3)$; ii) medium size data, $(m,n,r)\!=\!(50,1000,3)$; iii) large size data, $(m,n,r)\!=\!(100,10000,3)$.
Table \ref{tab:results_Exp1} reports the results on the datasets, according to the setting and metrics detailed before,
{where the highlighted rows correspond to the best-performing method and the proposed approach.
As detailed in the following analysis, RMU consistently achieves top-tier performance, ranking either first or second when jointly considering objective value, ranking, and sparsity.}

\begin{table}[h!]
\centering
\caption{The final function value, the ranking, the sparsity of H, obtained by the algorithms on 9 synthetic datasets.
}
\label{tab:results_Exp1}
\begin{tabular}{l|l|p{3cm} p{2.2cm} p{3.175cm}}
\hline
& \textbf{Method} & mean$\pm$std of $F_{\text{Last}}$ & ranking & $\H_{\text{Last}}$ median sparsity$\%$ 
\\
\hline
\multicolumn{5}{c}{small size data: \textbf{$(m,n,r,s,k_{\max},T_{\max})=(20,100,3,0.6,10^6,5)$}}
\\
\hline
\multirow{4}{*}{\rotatebox{90}{$\sigma=0$}} 
                  & EMUproj & 152.46$\pm$24.98& (13, 0, 8, 79)& 40.66\\
                  & SMUL1 & 150.30$\pm$27.67& (4, 1, 80, 15)& 73.16\\
                  & RMU & \textbf{141.52$\pm$25.39}& \textbf{(27, 56, 11, 6)}& \textbf{63.83}\\
                  & RCG & \textbf{140.04$\pm$25.42}& \textbf{(56, 43, 1, 0)}& \textbf{100.00}\\
\hline
\multirow{4}{*}{\rotatebox{90}{$\sigma=0.1$}} 
    & EMUproj & 172.91$\pm$25.44& (20, 0, 9, 71)& 41.00\\
    & SMUL1 &  177.29$\pm$27.48& (6, 0, 75, 19)& 83.33\\
    & RMU & \textbf{164.66$\pm$25.85}& \textbf{(53, 21, 16, 10)}& \textbf{62.66}\\
    & RCG & \textbf{163.27$\pm$26.07}& \textbf{(21, 79, 0, 0)}& \textbf{100.00}\\
\hline
\multirow{4}{*}{\rotatebox{90}{$\sigma=0.3$}} 
    & EMUproj & 235.00$\pm$28.44& (46, 0, 22, 32)& 39.33\\
    & SMUL1 & 257.26$\pm$34.77& (6, 0, 56, 38)& 88.00\\
    & RMU &  \textbf{232.73$\pm$31.29}& \textbf{(45, 3, 22, 30)}& \textbf{57.00}\\
    & RCG & \textbf{232.23$\pm$31.95}& \textbf{(3, 97, 0, 0)}& \textbf{100.00}\\
\hline
\multicolumn{5}{c}{medium size data: \textbf{$(m,n,r,s,k_{\max},T_{\max})=(50,1000,3,0.6,10^6,5)$}}
\\
\hline
\multirow{4}{*}{\rotatebox{90}{$\sigma=0$}} 
        & EMUproj &  \textbf{4997.57$\pm$371.04}& \textbf{(82, 0, 16, 2)}& \textbf{33.95}\\
        & SMUL1 & 5385.95$\pm$408.19& (5, 0, 81, 14)& 69.83\\
        & RMU & \textbf{5128.34$\pm$369.93}& \textbf{(10, 3, 3, 84)}& \textbf{59.11}\\
        & RCG & 5171.70$\pm$391.18& (3, 97, 0, 0)& 99.96\\
\hline
\multirow{4}{*}{\rotatebox{90}{$\sigma=0.1$}} 
        & EMUproj & \textbf{5644.77$\pm$424.27}& \textbf{(86, 0, 14, 0)}& \textbf{33.91}\\
        & SMUL1 &  6312.15$\pm$455.60& (4, 0, 86, 10)& 77.90\\
        & RMU & \textbf{5843.76$\pm$412.32}& \textbf{(10, 0, 0, 90)}& \textbf{55.76}\\
        & RCG & 5943.41$\pm$431.83& (0, 100, 0, 0)& 99.96\\
\hline
\multirow{4}{*}{\rotatebox{90}{$\sigma=0.3$}} 
        & EMUproj & \textbf{7696.90$\pm$433.29}& \textbf{(89, 0, 11, 0)}& \textbf{33.83}\\
        & SMUL1 & 9078.84$\pm$561.37& (8, 0, 89, 3)& 83.43\\
        & RMU &  \textbf{8018.14$\pm$499.01}& \textbf{(3, 0, 0, 97)}& \textbf{50.16}\\
        & RCG &  8327.87$\pm$540.22& (0, 100, 0, 0)& 100.00\\
\hline
\multicolumn{5}{c}{large size data: \textbf{$(m,n,r,s,k_{\max},T_{\max})=(100,10000,3,0.6,10^6,5)$}} \\
\hline
\multirow{4}{*}{\rotatebox{90}{$\sigma=0$}} 
    &EMUproj &  \textbf{133414.79$\pm$5835.10}& \textbf{(99, 0, 1, 0)}& \textbf{33.45}\\
    &SMUL1 & 153855.38$\pm$7467.90& (1, 0, 99, 0& 67.68\\
    &RMU & \textbf{142953.11$\pm$6565.51}& \textbf{(0, 99, 0, 1)}& \textbf{50.77}\\
    &RCG & 155081.24$\pm$8132.66& (0, 1, 0, 99)& 99.99\\
\hline
\multirow{4}{*}{\rotatebox{90}{$\sigma=0.1$}} 
    &EMUproj & \textbf{151988.78$\pm$6317.54} & \textbf{(98, 0, 2, 0)} &  \textbf{33.45}  \\
    &SMUL1 &  179816.74$\pm$8255.74& (2, 0, 98, 0)& 70.72\\
    &RMU & \textbf{161917.73$\pm$6982.22} & \textbf{(0, 98, 0, 2)} & \textbf{47.59} \\
    &RCG & 179762.51$\pm$8911.50& (0, 2, 0, 98)& 99.99\\
\hline
\multirow{4}{*}{\rotatebox{90}{$\sigma=0.3$}} 
    &EMUproj & \textbf{205028.43$\pm$8094.90}& \textbf{(95, 0, 5, 0)}& \textbf{33.47}\\
    &SMUL1 & 251804.87$\pm$10254.64& (5, 0, 95, 0)& 76.92\\
    &RMU &  \textbf{213833.68$\pm$8784.78}& \textbf{(0, 73, 0, 27)}& \textbf{40.87}\\
    &RCG &  238907.39$\pm$9286.03& (0, 27, 0, 73)& 99.99\\
\hline
\end{tabular}
\end{table}

Based on the ranking values shown in the tables, it can be observed that RCG performs better on small datasets, while RMU excels with larger datasets.
The proposed method is faster and more cost-efficient due to lower computational overhead per iteration.
Therefore, RMU can be considered a suitable alternative for large-scale problems when extremely high accuracy is not essential.

For the sparsity, we remark that, for RCG, the median sparsity\% is always 100\%.
This is because after say 39289 iterations, the smallest entry of $\H_{\text{Last}}$ produced by RCG has a magnitude of $5.3\!\times\!10^{-8}$, not reaching $10^{-9}$.
Regarding the comparison of other methods, we should remark that even when they reach a better sparsity measure or minimization point, these come from heuristic approaches where the normalization is performed a posteriori with a projection step.
In fact, although both methods preserve the normalization of the columns of $\H$, RMU achieves a better minimization of the objective function, while also yielding improved sparsity and integrating the normalization directly into the optimization process without requiring additional computational steps.

\begin{figure}[h!]
\centering\includegraphics[width=0.97\linewidth]{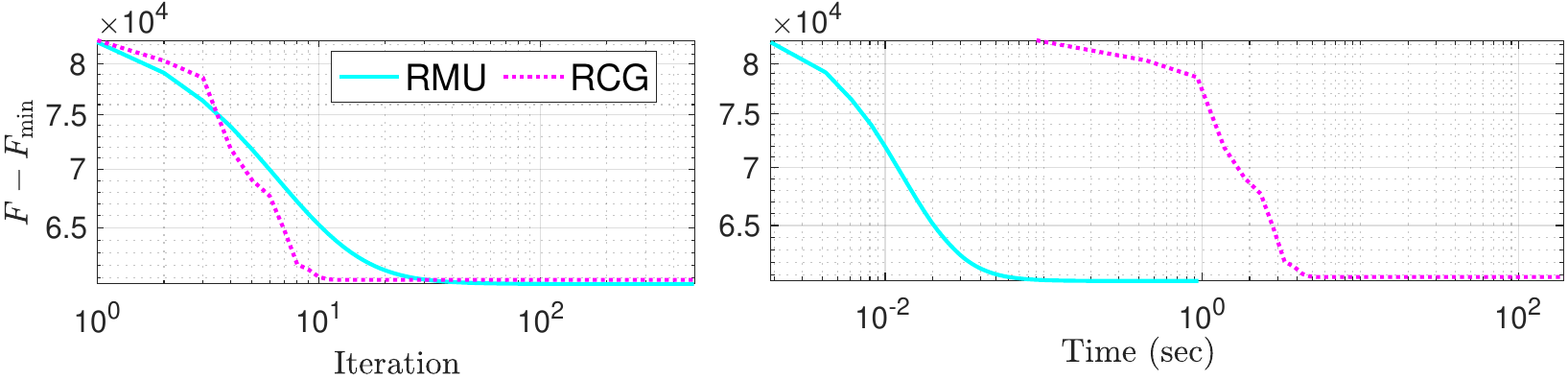}
\caption{Result on the hyperspectral dataset Samson with $m\!=\!156$ and $n\!=\!9025$ with $r\!=\!3$.
RMU takes about 10 seconds for 10000 iterations, while RCG takes 200 seconds for 500 iterations.
The iteration plot is showing the first 500 iterations.
}
\label{fig:Samson}
\end{figure}

\begin{figure}[h!]
\centering\includegraphics[width=0.75\linewidth]{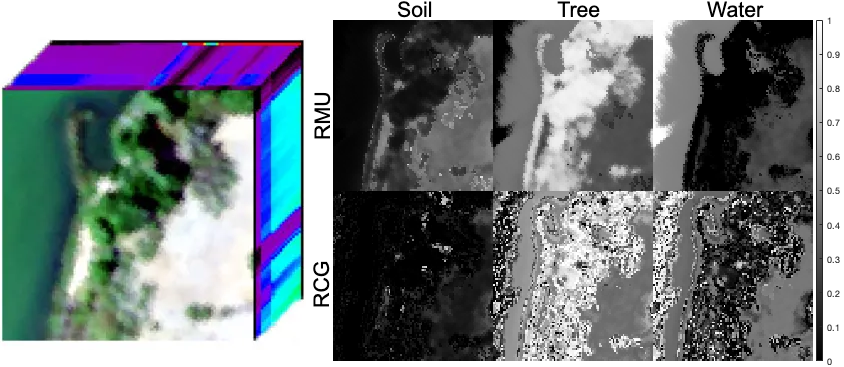}
\caption{Result on the abundance map on the decomposition of the hyperspectral dataset Samson.
Left: the Samson data cube.
Right: the decomposition.}
\label{fig:SamsonCube}
\end{figure}

\textbf{Simple preliminary experiment on hyperspectral dataset.}
We apply RMU and RCG on Samson hyperspectral dataset under the same setup condition.
We initialize $\H$ with NNSVD and we use the same penalization parameter $\lambda\!=\!0.35$.
As it can be seen in Fig.\ref{fig:Samson} in the Riemannian optimization framework RMU converges faster than RCG, moreover, Fig.\ref{fig:SamsonCube} shows that the proposed method is able to provide more interpretable decomposition of the three endmembers retained.
{ Specifically, the abundance maps generated by RMU present sharper spatial localization and reduced background noise compared to those typically produced by unconstrained approaches, making the physical distribution of the materials much clearer to identify.}

Although these are preliminary results, they already provide a good indication of the strong potential of RMU, suggesting that it may outperform existing approaches both in terms of speed and overall performance.
These findings motivate further investigation, and a detailed study of its practical applications is deferred to future work.

\section{Conclusion}\label{sec_conc}
In this work, we have proposed a novel manifold optimization method to solve low-rank problems with sparse simplex constraints by leveraging oblique manifolds.
Our approach ensures that the nonnegativity, sparsity, and sum-to-1 conditions are preserved throughout the optimization process while maintaining the convergence properties of Riemannian gradient descent.

Experiments on synthetic datasets have demonstrated the effectiveness of the proposed method compared to standard Euclidean and Riemannian optimization techniques.
Specifically, we obtained comparable, and sometimes superior, solutions faster than the baseline methods.
Moreover, unlike the RCG approach, the proposed RMU does not require computing vector transport or performing Wolfe line searches, making RMU appealing for large-scale problems.
In fact, our approach exhibits improved convergence behavior, numerical stability, and the ability to better exploit the underlying low-rank structure of the problem, particularly for high-dimensional problems.
These advantages make it particularly suitable for applications where structured low-rank solutions with simplex constraints are essential, such as in signal processing, machine learning, and computational biology.

Future work will focus on extending the method to explore adaptive strategies for manifold parameterization and investigating its applicability to real-world datasets in diverse domains.

\backmatter
\bmhead{Acknowledgements}
F.E. is member of the Gruppo Nazionale Calcolo
Scientifico - Istituto Nazionale di Alta Matematica (GNCS - INdAM).
F.E. thanks Prof. Del Buono from University of Bari for the useful scientific discussion.
The authors thank the Editor and the referee for their constructive comments and suggestions, which have improved the quality and presentation of this work.

\section*{Declarations}
\begin{itemize}
\item Funding: This work was supported by INdAM - GNCS Project \textit{Ottimizzazione e Algebra Lineare Randomizzata} (CUP E53C24001950001 to F.E.), and Piano Nazionale di Ripresa e Resilienza (PNRR), Missione 4 “Istruzione e Ricerca”-Componente C2 Investimento 1.1, “Fondo per il Programma Nazionale di Ricerca e Progetti di Rilevante Interesse Nazionale”, Progetto PRIN-2022 PNRR, P2022BLN38, \textit{Computational approaches for the integration of multi-omics data} (CUP H53D23008870001, to F.E.).
\item Conflict of interest/Competing interests (check journal-specific guidelines for which heading to use): The authors declare no competing interests.
\item Ethics approval and consent to participate: Not applicable.
\item Consent for publication: All authors give the consent for publication.
\item Data availability : data are synthetically generated.
\item Materials availability: Not applicable.
\item Code availability: The code will be accessible on \url{https://github.com/flaespo}.
\item Author contribution: All authors contributed equally to the work.
\end{itemize}
\bibliography{sn-bibliography}
\end{document}

%% file: MathSymbols.tex
\usepackage{bm}                 % enable bold-ed math symbol
\DeclareMathOperator*{\argmin}{\textrm{argmin}\,} % argmin
\def\diag            {\textrm{diag}}     % diag
\def\grad            {\textrm{grad}}     % grad   
\def\ones            {\bm{1}}            % All one vector or matrix
\def\proj            {\textrm{proj}}     % proj
\def\sgn             {\textrm{sign}}     % sign
\def\st              {\,\textrm{s.t.}\,} % such that (with spaces before/after)
\def\zeros           {\bm{0}}            % All zero vector or matrix
\def\IN {\mathbb{N}} 		% Natural number IN
\def\IR  {\mathbb{R}} 	    % Real set
\def\IRn {\IR^{n}} 			% Rn 
\def\IRmn{\IR^{m \times n}} % R m times n
\def\IRnn{\IR^{n \times n}} % R n-by-n matrix
\def\IRmr{\IR^{m \times r}} % R m-by-r matrix
\def\IRrn{\IR^{r \times n}} % R r-by-n matrix

\def\A {\bm{A}}
\def\B {\bm{B}}

\def\D {\bm{D}}
\def\E {\bm{E}}

\def\G {\bm{G}}
\def\H {\bm{H}} 
\def\I {\bm{I}}

\def\M {\bm{M}}

\def\P {\bm{P}}
\def\Q {\bm{Q}}

\def\V {\bm{V}}
\def\W {\bm{W}}
\def\X {\bm{X}} 

\def\Z {\bm{Z}}

\def\x {\bm{x}}

\def\cB { \mathcal{B}}

\def\cM { \mathcal{M}}
\def\cN { \mathcal{N}}
\def\cO { \mathcal{O}}
\def\cP { \mathcal{P}}

\def\cR { \mathcal{R}}
\def\cS { \mathcal{S}}
\def\cT { \mathcal{T}}

\def\balpha   {\boldsymbol{\alpha}}